\documentclass[conference, 10pt]{IEEEtran}
\usepackage[english]{babel}
\usepackage{graphicx}
\usepackage{amsfonts,amssymb,amsmath}
\usepackage{hyperref} 
\usepackage{comment} 
\usepackage[dvipsnames]{xcolor} 
\usepackage{soul}  
\usepackage{float} 
\usepackage{cite} 
\usepackage{overpic} 
\usepackage{todonotes}
\usepackage{siunitx}
\usepackage{ulem} 

\usepackage{ifthen}

\newcommand{\diff}{\,\mathrm{d}}

\newcommand{\nh}{\vv{\hat{n}}}

\newcommand{\Tc}{ { \mathcal{T}} }

\newcommand{\Kc}{ { \mathcal{K}} }

\newcommand{\Rc}{ { \mathcal{R}} }

\newcommand{\vv}[1]{\ensuremath{\boldsymbol{#1}}}           
\newcommand{\vet}[1]{\mathbf{#1}}
\def\Xtil{\vet{{\Tilde{X}}}}
\def\Btil{\vet{{\Tilde{B}}}}

\input{Authors_and_title}

\begin{document}
\newmaketitle

\begin{abstract}
Inverse source approaches have shown their relevance for several applications in the past years. They rely on the solution of an ill-posed problem where near-field/current data is reconstructed starting from far-field (or less informative field) information. Standard strategies, including the physically constrained ones using Love conditions, result in linear systems to be pseudoinverted which are still ill-conditioned due to the lack of information from the evanescent fields. In this work we present a generalized pseudoinverse for those problems that allows the inclusion of extra constraints from the evanescent field space, when available. This is obtained by dropping some of the standard Moore-Penrose (MP) requirements and using the resulting degrees of freedom to obtain a generalized pseudoinverse that shows favorable performance in several cases of practical interest.
\end{abstract}

\section{Introduction}\label{sec:intro}
Electromagnetic inverse source problems have been effectively used for tackling  a variety of application scenarios including antenna diagnostics and far to near field reconstructions (see \cite{Vecchi} and references therein). In light of the equivalence theorem, many formulations are allowed and proved to be comparable, in reconstructing the far-field at least, in terms of speed and accuracy once discretized and solved \cite{Korn}.
However, all methods inherit the nature of the same electromagnetic problem and describe a common ill-posedness which gets more severe as the distance from the source(s) increases. In fact, when dealing with far-field measurements, the inevitable loss of information associated with the evanescent modes compromises near-field reconstruction: all the measured power concentrates in the well discretized propagating modes \cite{Hansen}. On the one hand this problem has been partially mitigated using algebraic pseudoinverses, most notably the Moore-Penrose pseudoinverse (MPP), associated with constraints like the Love one \cite{Vecchi}, which however has limited impact on the near-field reconstruction given that the regularization is done with information from the internal equivalent problem. Other approaches, often in the context of near to far field transformations and pattern characterization (see the nice work \cite{Migliore} and related papers) have been instead focusing in recovering overcomplete solution spaces in which the reconstruction is performed.
In this work, we propose a different approach which is somehow in between the two. A limited number of physical constraints (a substantially undercomplete basis) is considered and used to define a generalized pseudoinverse that still allows for norm minimization or equivalent conditions and, most importantly, is still compatible with Love constraints that are so popular in the community. The pseudoinverse we propose compensates for the loss of evanescent degrees of freedom with a limited number of constraints obtained by measurements or simulation while still solving the better conditioned part of the matrix through an algebraic pseudoinversion. The new pseudoinverse is compatible with most of the existing inverse source formulations (where it can replace the MPP or equivalent). Theoretical considerations are complemented by numerical results showing the practical relevance of the proposed approach.

\section{Notation and Background}\label{sec:notation}
A given source can be enclosed in a fictitious equivalent surface $\Gamma$ which delimits an internal/external volume $\Omega^-/\Omega^+$. While $\Omega^+$ retains the properties of the original medium, $\Omega^-$ can be filled with whatever material; here the same permittivity $\epsilon_0$ and permeability $\mu_0$ as $\Omega^+$ are chosen. The equivalence theorem allows us to find a set of surface currents densities $\vv{M} = \nh \times\left(\vv{E^-}-\vv{E^+}\right)$ and $\vv{J} = \nh\times\left(\vv{H^+}-\vv{H^-}\right)$ on $\Gamma$ that radiate the same fields as the source, where $\nh$ is the unit normal vector field on $\Gamma$ pointing into $\Omega^+$ and $\vv{E^+},\,\vv{H^+}$ (resp. $\vv{E}^-,\,\vv{H}^-$) are the electric and magnetic field in $\Omega^+$ (resp. $\Omega^-$).

We then define the electric operator on $\Gamma$ as $\Tc(\vv{J})(\vv{r})=ik\,\Tc_s(\vv{J})(\vv{r})-(ik)^{-1}\,\Tc_h(\vv{J})(\vv{r})$ with $\vv{r}\in\Gamma$, $k$ the wavenumber,
$\Tc_s(\vv{J})(\vv{r}) = \nh\times\int_{\Gamma}\frac{e^{ik\left|\vv{r}-\vv{r}'\right|}}{4\pi\left|\vv{r}-\vv{r}'\right|}{\vv{J}(\vv{r}')}\diff\vv{r}'$, and $\Tc_h(\vv{J})(\vv{r}) = \nh\times\nabla\int_{\Gamma}\frac{e^{ik\left|\vv{r}-\vv{r}'\right|}}{4\pi\left|\vv{r}-\vv{r}'\right|}{\nabla_s\cdot\vv{J}(\vv{r}')}\diff\vv{r}'$.
The magnetic operator is $\Kc(\vv{J})(\vv{r})  =  - \nh \times p.v. \int_\Gamma
    \nabla \times \frac{e^{i k \left|\vv{r}-\vv{r}'\right|}}{4 \pi \left|\vv{r}-\vv{r}'\right|} \vv{J}(\vv{r}')\diff\vv{r}'$
with $\vv{r}\in \Gamma$. We extend these definition to the case in which $\vv r$ belongs to a measurement surface $\Gamma_m$ distinct from $\Gamma$ and denote these extensions as $\Tc_r$ and $\Kc_r$; their expression are similar to $\Tc$ and $\Kc$ except that the principal value in $\Kc$ is not necessary. Finally, given a measurement vector $b\in\Gamma_m$, the tangentially projected radiation operator is
\begin{equation}\label{eq:Rdef}
    \Rc=\begin{bmatrix}
        -\Kc_r & \Tc_r\\
        -\Tc_r &-\Kc_r
    \end{bmatrix}\,,
    \quad
    \vv{r}\in\Gamma_m\,,
\end{equation}
and its related linear system is
\begin{equation}\label{eq:ls}
    \Rc x=\Rc\cdot\begin{bmatrix}-\vv{M}\\\eta_0\vv{J}\end{bmatrix}=b=\nh\times\begin{bmatrix}\vv{E}^+\\\eta_0\vv{H}^+\end{bmatrix}
\end{equation}
with $\eta_0=\sqrt{\mu_0/\epsilon_0}$. After discretization \eqref{eq:ls} reads
\begin{equation}\label{eq:ls_disc}
    \vet{Rx} = \vet{b}
\end{equation}
with $\vet{R}$, $\vet{x}$ and $\vet{b}$ being the Rao-Wilton-Glisson discretized version of $R$, $x$ and $b$ respectively.

\section{A Constrained Pseudoinverse}\label{sec:fpa_pinv}

Solving the inverse source problem amounts to pseudoinverting $\vet{R}$ to solve the ill-posed problem \eqref{eq:ls_disc}. This can be done after adding direct or indirect constraints, including single source enforced solutions, Love, and/or minimum norm conditions (see \cite{Korn,Vecchi} and references therein). Often, when $\vet{R}$ (or its regularized version) is pseudoinverted, the MPP or an equivalent iterative pseudoinverse is used. This norm minimizing approach warranties four fundamental properties for the pseudoinverse, i.e. the MPP pseudoinverse $\vet{R}^+$ of $\vet{R}$ satisfies 1) $\vet{R}\vet{R}^+\vet{R}=\vet{R}$, 2) $\vet{R}^+\vet{R}\vet{R}^+=\vet{R}^+$, 3) $(\vet{R}^+\vet{R})^*=\vet{R}^+\vet{R}$, and 4) $(\vet{R}\vet{R}^+)^*=\vet{R}\vet{R}^+$. The MPP is the unique pseudoinverse that satisfies 1)--4), however, among the above properties only 1) is unavoidable if we want to find a solution of \eqref{eq:ls_disc}; thus renouncing to some of the properties 2)-4) may allow to further optimize the pseudoinverse. A natural starting point is the general pseudoinverse $\vet{R}^\dagger$ of $\vet{R}$ that only satisfies 1). This would provide an exact solution for \eqref{eq:ls_disc} that would allow the injection of \textit{a priori} regularity constraints and minimization criteria,
but the selection criterion of the solution may not be sufficient to describe the physics of the problem. Another pseudoinverse that could be of interest would be the one providing an approximation of the solution of \eqref{eq:ls_disc}, for a specific set of linearly independent right hand sides $\Btil$, in the space of linear combinations of the corresponding known solutions $\Xtil$. 
One way of designing such a pseudoinverse is $\Xtil\Btil^{\bullet}$ where $\Btil^{\bullet}$ is a pseudoinverse of $\Btil$ that satisfies 1). 
Although $\Xtil\Btil^{\bullet}$ provides physically relevant solutions, when the number of vectors is lower than the degrees of freedom, it provides approximate solutions that are not necessarily solving \eqref{eq:ls_disc}. In fact, if $\Btil^{\bullet}$ is a pseudoinverse satisfying at least 1) then $\Xtil\Btil^{\bullet}$ is a pseudoinverse satisfying at least 2), so it is not ensured to yield an exact solution of \eqref{eq:ls_disc}.
In this work we propose a suitable operator average of the two approaches, i.e. a more general pseudoinverse in between the two described above:
\begin{equation}\label{eq:int_pinv}
    \vet{R}^\ddag = \vet{R}^\dag(\vet{I}-\vet{\Btil\Btil}^\bullet) + \vet{\Xtil\Btil}^\bullet.
\end{equation}
The above pseudoinverse has the advantage of providing a part of the solution as a linear combination of know solutions vectors but, at the same time and differently from $\Xtil\Btil^{\bullet}$, provides an exact reconstruction of the measurement data since it satisfies 1) and, in particular, $\vet{R}(\vet{R}^\ddag\vet{b}) = \vet{R}(\vet{R}^\ddag\vet{R}\vet{x}) = \vet{b}$.

\section{Numerical Results}\label{sec:numerical}
\begin{figure}
    \centering
    \includegraphics[height=5cm]{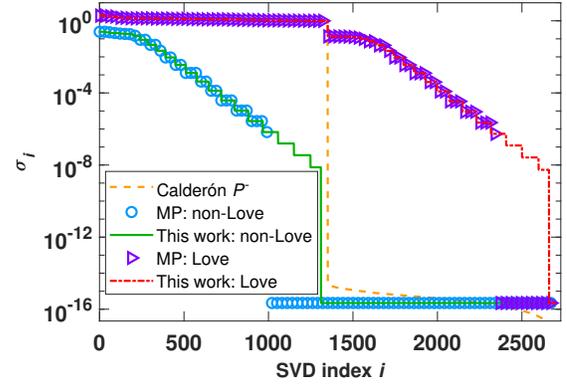}
    \caption{Spectral analysis of the different formulations.}\label{fig:rank_enhance}
\vspace{-0.15cm} 
\end{figure}

\begin{figure}
    \centering
    \begin{overpic}[height=5cm,percent]{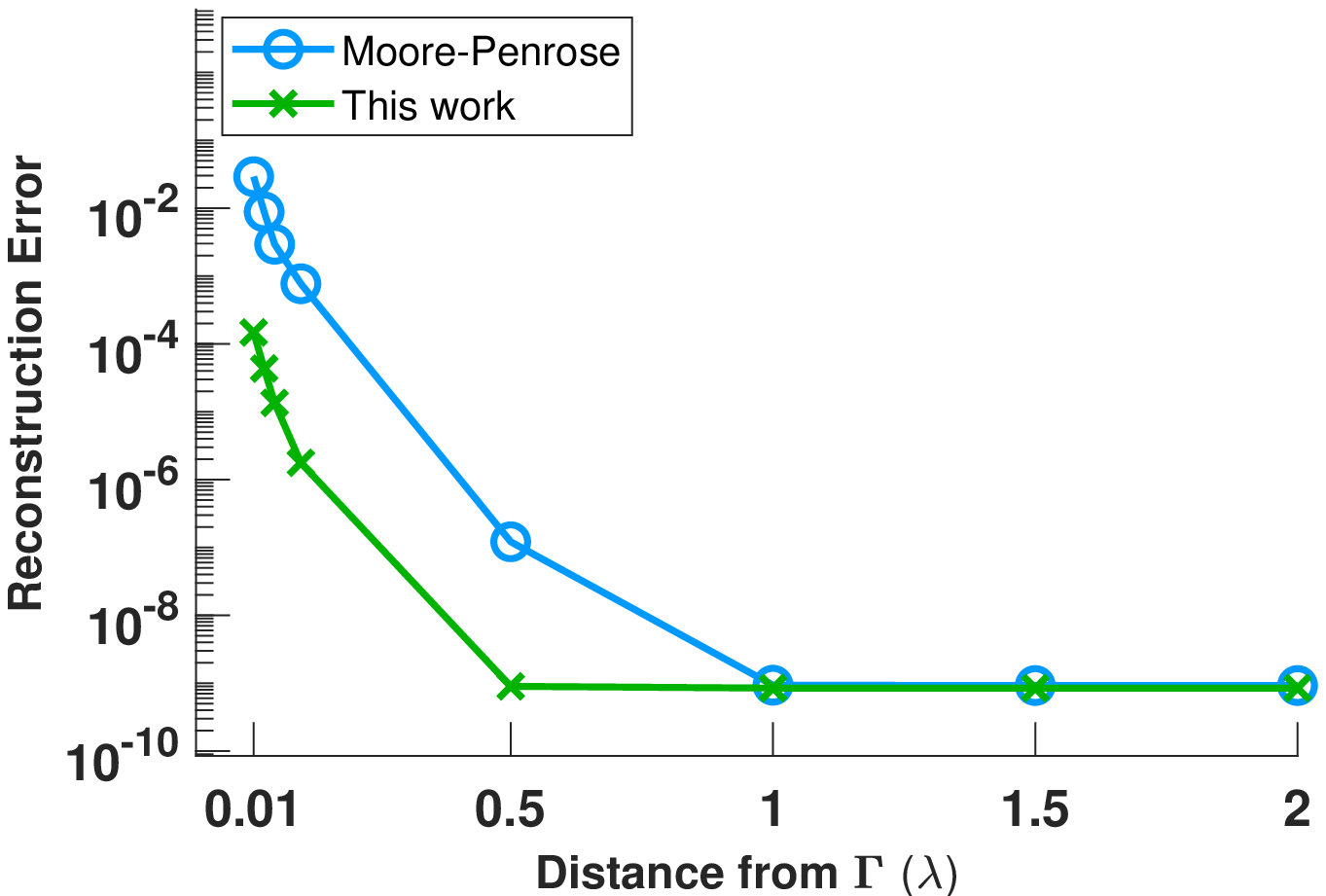}
    \put(35,18){
    \includegraphics[scale=0.22]{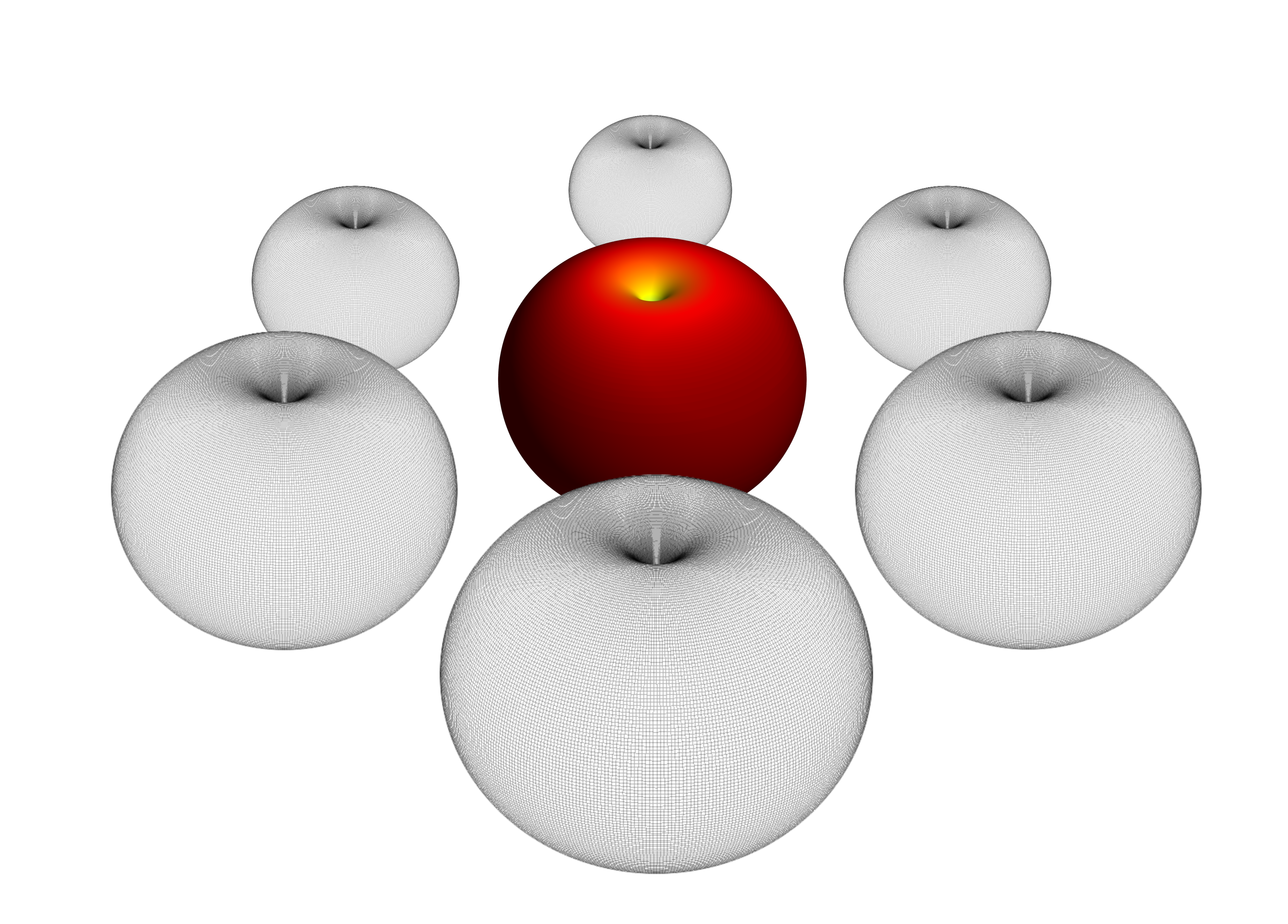}}
    \end{overpic}
    \caption{Original (red) and constraining (grey) dipole sources in 3D space and associated reconstruction error from far-field synthetic measurements.}\label{fig:loc_interp}
\vspace{-0.35cm} 
\end{figure}

The constrained pseudoinverse we propose has been numerically analysed. First the rank-enhancement effect of $\mathbf{R}^\ddag$ with respect to $\mathbf{R}^+$ is shown in Fig. \ref{fig:rank_enhance} through a semi-analytic spectral analysis realized with vector spherical harmonics that shows that the novel pseudoinverse is able to recover the evanescent base truncated by the MPP.
To highlight that the regularization is near-field aware (differently from other forms of regularization) we have studied a Hertzian dipole oscillating at \SI{5}{\giga\hertz} with a measurement surface at a distance of two wavelengths from the reconstruction sphere of radius \SI{4}{\centi\meter} and we have discretized \eqref{eq:ls_disc} using the Galerkin method. The measurement vector is pseudoinverted both with the MPP and the constrained pseudoinverse through a near-field spherical placement of identical sources. A simplified representation of the spatial constraining along with reconstruction results is shown in Fig. \ref{fig:loc_interp}. Compared to the MPP, $\vet{R}^\ddag$ shows an improvement in near-field reconstruction while not altering the far-field accuracy, in line with our theoretical considerations.

\section*{ACKNOWLEDGMENT}\label{sec:ack}
This work was supported in part by the European Research Council (ERC) through the European Union’s Horizon 2020 Research and Innovation Programme under Grant 724846 (Project 321); in part by the Italian Ministry of University and Research within the Program PRIN2017, EMVISIONING, under Grant 2017HZJXSZ, CUP:E64I190025300; in part by the Italian Ministry of University and Research within the Program FARE, CELER, under Grant R187PMFXA4.
\vspace{-0.05cm} 

\begin{thebibliography}{1}

\bibitem{Korn} J. Kornprobst, R. A. M. Mauermayer, O. Neitz, J. Knapp, and T. F. Eibert, ``On the Solution of Inverse Equivalent Surface-Source Problems'', PIER, vol. 165, pp. 47–65, 2019.
\bibitem{Vecchi} J. L. A. Quijano and G. Vecchi, ``Field and Source Equivalence in Source Reconstruction on 3D Surfaces'', PIER, vol. 103, pp. 67–100, 2010.
\bibitem{Hansen} J. E. Hansen, ``Spherical Near-ﬁeld Antenna Measurements'',
Vol. 26, IEE Electromagnetic Waves Series, Stevenage Herts
England Peter Peregrinus Ltd., pp. 8-26, 1988.
\bibitem{Migliore} M. Salucci, M. D. Migliore, G. Oliveri, and A. Massa, ``Antenna Measurements-by-Design for Antenna Qualification'', IEEE Trans. Antennas Propagat., vol. 66, no. 11, pp. 6300–6312, Nov. 2018.


\end{thebibliography}

\end{document}